\documentclass{article}  

\usepackage{amsmath, amssymb, amsfonts, amsbsy, latexsym}
\usepackage{graphicx,stmaryrd}

\newtheorem{theorem}{Theorem}[section]

\newtheorem{lemma}[theorem]{Lemma}
\newtheorem{proposition}[theorem]{Proposition}

\newtheorem{remark}[theorem]{Remark}
\newtheorem{example}[theorem]{Example}

\def\span{{\mbox{\rm span}}}

\def\C{\mathbb{C}}
\def\R{\mathbb{R}}

\def\QQ{\mathbb{Q}}

\def\<{\langle}
\def\>{\rangle}

\def\epsilon{\varepsilon}

\newcommand{\ip}[2]{\langle#1,#2\rangle}

\newcommand{\norm}[1]{\left\lVert#1\right\rVert}

\newcommand{\fc}{{\mathcal F}}
\newcommand{\gc}{{\mathcal G}}
\newcommand{\hc}{{\mathcal H}}

\newcommand{\nonlin}{\alpha}
\newcommand{\phik}{{\varphi_k}}
\newcommand{\outp}[2]{{\llbracket #1,#2 \rrbracket }}
\newcommand{\ac}{\mathcal{A}}
\renewcommand{\SS}{\mathcal{S}}
\newcommand{\Sonezero}{{\SS^{1,0}}}
\newcommand{\Soneone}{{\SS^{1,1}}}
\newcommand{\xb}{{\mathbf x}}
\renewcommand{\j}{{\mathbf j}}

\title{Stability of Phase Retrievable Frames}

\author{Radu Balan \\
Mathematics Department and CSCAMM, University of Maryland,\\
College Park, MD 20742, USA
}

\begin{document}
  \maketitle

\begin{abstract}
In this paper we study the property of phase retrievability by redundant sysems of vectors under perturbations of the frame set.
Specifically we show that if a set $\fc$ of $m$ vectors in the complex Hilbert space of dimension n allows for vector reconstruction from magnitudes of
its coefficients, then there is a perturbation bound $\rho$ so that any frame set within $\rho$ from $\fc$ has the same property. In particular
this proves the recent construction in \cite{BH13} is stable under perturbations. By the same token we reduce the critical
cardinality conjectured in \cite{BCMN13a} to proving a stability result for non phase-retrievable frames.
\end{abstract}



\section{INTRODUCTION}
\label{sec:intro}

The {\em phase retrieval} problem presents itself in many applications is physics and engineering. Recent papers on this topic
(see \cite{BCE06,CSV12,Bal12a,Bal13a,ABFM12,BCMN13a,W12}) present a full list of examples ranging from X-Ray crystallography to
 audio and image signal processing, classification with deep networks, quantum information theory, and fiber optics data transmission.

In this paper we consider the complex case, namely the Hilbert space $H=\C^n$ endowed with the usual Euclidian scalar product
 $\ip{x}{y}=\sum_{k=1}^n x_k\overline{y_k}$.
On $H$ we consider the equivalence relation $\sim$ between two vectors $x,y\in H$ defined as follows; the vectors $x$ and $y$ are similar
$x\sim y$ if and only if there is a complex constant $z$ of unit magnitude, $|z|=1$, so that $y=zx$. Let $\hat{H}=H/\sim$ be the
quotient space. Thus an equivalence class (a {\em ray}) has the form $\hat{x}=\{e^{i\varphi}x~,~\varphi\in[0,2\pi) \}$. 
A subset $\fc\subset H$ of the Hilbert space $H$ (regardless whether it is finite dimensional or not) is called {\em frame} if there are two positive 
constants $0<A\leq B<\infty$ (called {\em frame bounds}) so that for any vector $x\in H$, 
\begin{equation}
\label{eq:frame}
 A\norm{x}^2 \leq \sum_{f\in \fc} |\ip{x}{f}|^2 \leq B \norm{x}^2 
\end{equation}
In the finite dimensional case considered in this paper, the frame condition simply reduces to the spanning condition.  
Specifically $\fc=\{f_1,\ldots,f_m\}$ is frame for $H$ if and only if $H=span(\fc)$. Obviously $m\geq n$ must hold. 
When we can choose $A=B$ the frame is called {\em tight}. If furthermore $A=B=1$ then $\fc$ is said a {\em Parseval frame}.
Consider the following nonlinear map
\begin{equation}
\label{eq:alpha}
\nonlin : \hat{H}\rightarrow \R^m~~,~~(\nonlin(\hat{x}))_k=|\ip{x}{f_k}|~,~1\leq k\leq m
\end{equation}
which is well defined on the classes $\hat{x}$ since $|\ip{x}{f_k}|=|\ip{y}{f_k}|$ when $x\sim y$.

The frame $\fc=\{f_1,\ldots,f_m\}$ is called {\em phase retrievable} is the nonlinear map $\nonlin$ is injective.
Notice that any signal $x\in H$ is uniquely defined by the magnitudes of its frame coefficients $\nonlin(x)$ up to a global phase factor,
if and only if $\fc$ is phase retrievable. The main result of this paper states that the phase retrievable property is stable under small
perturbations of the frame set. Specifically we show

\begin{theorem}
\label{th1}
Assume $\fc=\{f_1,\ldots,f_m\}$ is a phase retrievable frame for a complex Hilbert space $H$. Then there is a $\rho>0$ so that
any set $\fc'=\{f_1',\ldots,f_m'\}$ with $\norm{f_k-f_k'}<\rho$, $1\leq k\leq m$, is also a phase retrievable frame.
\end{theorem}

We prove this theorem in section \ref{sec2}. The proof is based on a recent necessary and sufficent condition obtained independently
in \cite{BCMN13a} and \cite{Bal13a}. The exact form of this result is slightly different than the equivalent results stated in the aforementioned papers.
Consequently we will provide a direct proof. 

An interesting problem on phase retrievable frames is to find a critical cardinal $m^*(n)$ that has the following properties:
\begin{itemize}
\item[(A)] For any $m\geq m^*(n)$ the set of phase retrievable frames is generic with respect to the Zariski topology;
\item[(B)] If $\fc$ is a phase retrievable frame of $m$ vectors, then $m\geq m^*(n)$.
\end{itemize}
Clearly (B) is equivalent to:
\begin{itemize}
\item[(C)] If $m<m^*(n)$ there is no frame $\fc$ of $m$ vectors that is also phase retrievable.
\end{itemize}
The current state-of-the-art on this problem is summarized by the
following statements:
\begin{enumerate}
\item see [\cite{BCE06}]. If $m\geq 4n-2$ then generically (with respect to the Zariski topology) any frame is phase retrievable for $\C^n$;
\item see [\cite{MV13}]. For generic 4 $n\times n$ unitary matrices on $\C^n$, any subset of $m=4n-3$ columns forms a phase retrievable frame;
\item see [\cite{HMW11}]. If $\fc$ is a phase retrievable frame in $\C^n$ then
\[ m\geq 4n-2-2\beta + \left\{ \begin{array}{rl}
2 & \mbox{$if~n~odd~and~\beta=3\,mod\,4$} \\
1 & \mbox{$if~n~odd~and~\beta=2\,mod\,4$} \\
0 & otherwise
\end{array} \right.
\]
where $\beta=\beta(n)$ is the number of $1$'s in the binary expansion of $n-1$.
\end{enumerate}
Hence, if such a critical cardinal exists, we know $4n-O(log(n))\leq m^*_{\C}(n) \leq 4n-2$.
The authors of \cite{BCMN13a} conjectured that $m^*_{\C}(n)=4n-4$. 
In the case $m=4n-4$, Bodmann and Hammen constructed \cite{BH13} a phase retrievable frame. In section \ref{sec3} we review their construction 
and we show it is stable under small perturbations. In section \ref{sec4} we consider the critical cardinality conjecture
and show that (C) is equivalent to a stability result for sets of frames that fail to be phase retrievable.

Note the corresponding problems for the real case are completely solved. In fact \cite{BCE06} gives a geometric condition equivalent to a frame being
phase retrievable in $\R^n$. That condition (namely, for any partition of the frame set $\fc=\fc_1\cup\fc_2$, at least one of $\fc_1$ or $\fc_2$ 
must span $\R^n$) is stable under small perturbations. Additionally that same condition implies that $m^*_{\R}=2n-1$ is the critical cardinal
in the real case.

\section{Notations}
\label{sec1b}

In this section we recall some notations we introduced in \cite{Bal13a} that will be used in the following sections. 
Let $\fc=\{f_1,\ldots,f_m\}$ be a frame in $H=\C^n$. Let $\j:H\rightarrow\R^{2n}$ denote the embedding 
\begin{equation} \j(x)=\left[ \begin{array}{c}
real(x) \\
imag(x)
\end{array} \right] 
\end{equation}
which is a unitary isomorphism $\j$ between $H$ seen as a real vector space endowed with the real inner product $\ip{x}{y}_{\R}=real(\ip{x}{y})$ and $\R^{2n}$:
\begin{equation} \ip{x}{y}_{\R} = real(\ip{x}{y}) = \ip{\j(x)}{\j(y)}. \end{equation}

For two vectors $u,v\in\R^{2n}$, $\outp{u}{v}$ denotes the symmetric outer poduct 
\begin{equation} \outp{u}{v}=\frac{1}{2}(uv^T+vu^T).  \end{equation}
and similarly for two vector $x,y\in \C^n$ denote by $\outp{x}{y}$ their symmetric outer product defined by
\begin{equation}
\outp{x}{y} = \frac{1}{2}(xy^*+yx^*)
\end{equation}
For each $n$-vector $f_k$ we denote by $\varphi_k$ the $2n$ real vector,
 and by $\Phi_k$ the symmetric nonnegative rank-2, $2n\times 2n$ matrix defined by
\begin{equation}
\label{phis}
\phik = \j(f_k) =\left[ \begin{array}{c}
\mbox{real($f_k$)} \\
\mbox{imag($f_k$)} 
\end{array} \right] ~~,~~
\Phi_k = \outp{\phik}{\phik}+\outp{J\phik}{J\phik}=\phik \phik^T + J\phik \phik^T J^T
~~,~~{\rm where}~J = \left[ 
\begin{array}{cc}
0 & - I  \\
I & 0 
\end{array} \right]
\end{equation}
Note the following key relations:
\begin{eqnarray}
real(\ip{x}{f_k}) & = & \ip{\xi}{\phik} \\
|\ip{x}{f_k}|^2 & = & \ip{\Phi_k \xi}{\xi}\\
real(\ip{x}{f_k}\ip{f_k}{y}) & = & \ip{\Phi_k \xi}{\eta}
\end{eqnarray}
where $\xi=\j(x)$ and $\eta=\j(y)$. For every $\xi\in\R^{2n}$ set
\begin{equation}
R(\xi) = \sum_{k=1}^m \outp{\Phi_k \xi}{\Phi_k\xi}
\end{equation}

Let $\Sonezero$ and $\Soneone$ denote the following spaces of symmetric operators over a Hilbert space $K$ (real or complex)
\begin{eqnarray}
\Sonezero(K) & = & \{ T\in Sym(K)~~,~~rank(T)\leq 1~,~\lambda_{max}(T)\geq 0 = \lambda_{min}(T) \} \\
\Soneone(K) & = & \{ T\in Sym(K)~~,~~rank(T)\leq 2 ~,~Sp(T)=\{\lambda_{max},0,\lambda_{min}\}~,~\lambda_{max}\geq 0
\geq \lambda_{min} \}
\end{eqnarray}
where $Sym(K)$ denotes the set of symmetric operators (matrices) on $K$, $Sp(T)$ denotes the spectrum (set of eigenvalues) of $T$,
 and $\lambda_{max},\lambda_{min}$ denote the largest, and smallest eigenvalue of the corresponding operator. Note
\[ \Sonezero(K) = \{ T=\outp{x}{x} ~~,~~x\in K \} \]
For the frame $\fc=\{f_1,\ldots,f_m\}$ we let $\ac$ denote the linear operator
\begin{equation}
\ac: Sym(H)\rightarrow \R^m~~,~~(\ac(T))_k = \ip{Tf_k}{f_k}=trace\left(T\outp{f_k}{f_k}\right)
\end{equation}
Note the frame condition (\ref{eq:frame}) reads equivalently:
\begin{equation}
\label{eq:Ac}
 A \norm{\outp{x}{x}}_1 \leq \norm{\ac(\outp{x}{x})}^2 \leq B \norm{\outp{x}{x}}_1 
\end{equation}
where $\norm{T}_1=\sum_{k=1}^{rank(T)}|\lambda_k(T)|$ denotes the nuclear norm of operator $T$, that is the sum of its
singular values, or the sum of magnitudes of its eigenalues when $T$ is symmetric.
The upper bound is obviously always true (for an appropriate $B$) in the case of finite frames. The lower bound
in (\ref{eq:frame}) or (\ref{eq:Ac}) is equivalent to the spanning condition  $\span(\fc)=H$.
In turn this spanning condition can be restated in terms of a null space condition for $\ac$. Specifically let $\ker(\ac)=\{T\in Sym(H)~,~\ac(T)=0\}$
denote the kernel of the linear operator $\ac$. Then $\span(\fc)=H$ (and therefore $\fc$ is frame) if and only if 
\begin{equation}
\label{eq:ker1}
\ker(\ac) \cap \Sonezero = \{ 0 \}
\end{equation}
Recall the nonlinear map $\nonlin$ introduced in (\ref{eq:alpha}). We shall consider also the square map $\nonlin^2$ defined by:
\begin{equation}
\label{eq:alpha2}
\nonlin^2: H\rightarrow\R^m~~,~~(\nonlin^2(x))_k = |\ip{x}{f_k}|^2
\end{equation}
Of course $\nonlin$ is injective if and only $\nonlin^2$ is injective.

\section{Stability of Phase Retrievable Frames}
\label{sec2}
We start by presenting two lemmas regarding the objects we introduced earlier.
\begin{lemma}
\label{lem1}
Fix a Hilbert space $K$.
\begin{enumerate}
\item As sets, $\Soneone(K) = \Sonezero(K) - \Sonezero(K)$. Explicitely this means:
\[ \forall\, T\in\Soneone ~\exists\, T_1,T_2\in\Sonezero~~s.t.~T=T_1-T_2~~~;~~~{\rm Conversely}:~\forall\, T_1,T_2\in\Sonezero~,~T_1-T_2\in\Soneone; \]
\item For any $T\in\Soneone(K)$ there are $u,v\in K$ so that $T=\outp{u}{v}$;
\item For any $u,v\in K$, $\outp{u}{v}\in\Soneone(K)$;
\item $\Soneone(K) = \{ T=\outp{u}{v}~,~u,v\in K\}$.
\end{enumerate}
\end{lemma}

The proof of this lemma can be found in \cite{Bal13a} Lemma 3.7.

\begin{lemma}
\label{lem2}
The following are equivalent:
\begin{enumerate}
\item The nonlinear map $\nonlin$ is injective;
\item $\ker(\ac)  \cap \Soneone = \{0 \}$
\item There is a constant $a_0>0$ so that 
\begin{equation}
\sum_{k=1}^m \left| |\ip{x}{f_k}|^2 - |\ip{y}{f_k}|^2 \right|^2 \geq a_0 \left( \norm{x-y}^2\norm{x+y}^2 + 4(imag(\ip{x}{y}))^2 \right)
\end{equation}
for any $x,y\in H\in\C^n$;
\item There is a constant $a_0>0$ so that for all $\xi\in\R^{2n}$, $\lambda_{2n-1}(R(\xi))\geq a_0 \norm{\xi}^2$  (here, $\lambda_{2n-1}(T)$ denotes the $2n-1^{th}$ largest
eigenvalue of $T$);
\item There is a constant $a_1>0$ so that for all $\xi\in\R^{2n}$,
\begin{equation}
\label{eq:R2}
L(\xi):=R(\xi) + J\outp{\xi}{\xi}J^T = \sum_{k=1}^m \Phi_k \xi\xi^T \Phi_k + J\xi\xi^TJ^T \geq a_1 \norm{\xi}^2 1_{\R^{2n}}
\end{equation} 
where the inequality is in the sense of quadratic forms;
\item For every $\xi\in\R^{2n}$, $\dim\ker R(\xi)=1$ ;
\item For every $\xi\in\R^{2n}$, $rank(R(\xi))=2n-1$ ;
\item For every $\xi\in\R^{2n-1}$, $\ker(R(\xi)) = \{ c\, J\xi~,~c\in\R\}$;
\item There is a constant $a_0>0$ so that for all $\xi\in\R^{2n}$,
\begin{equation}
\label{eq:proj}
R(\xi) \geq a_0 \norm{\xi}^2 (1-P_{J\xi})
\end{equation}
where $P_{J\xi} = \frac{1}{\norm{\xi}^2} \outp{J\xi}{J\xi}$ is the orthogonal projection onto $J\xi$.
\end{enumerate}
\end{lemma}
\begin{remark}
Before presenting the proof, note the constants $a_0$ at (iii), (iv) and (ix) can be chosen to be equal. Hence the optimal (i.e. the largest) $a_0$ is given by
\begin{equation}
\label{eq:a0}
a_0 = \min_{\norm{\xi}=1} \lambda_{2n-1}(R(\xi))
\end{equation}
Additionally, the constant $a_1$ at (v) can be chosen as $a_1=min(1,a_0)$.
\end{remark}

{\bf Proof of Lemma \ref{lem2} }

Claims (i),(ii),(iv),(vi)-(ix) are known in literature - see for instance \cite{Bal13a} Theorem 2.2 and the bibliographical indications - and Theorem 3.1 of \cite{Bal13a}.
Claim (v) follows from (ix) by adding $\norm{\xi}^2P_{J\xi}$ on both sides. Claim (iii) follows from Theorem 3.1 (2) of \cite{Bal13a}, where we set $u=x-y$ and $v=x+y$ and
by remarking $imag(\ip{u}{v}) = 2\,imag(\ip{x}{y})$ and $real(\ip{u}{f_k}\ip{f_k}{v}) = real(|\ip{x}{f_k}|^2-|\ip{y}{f_k}|^2)$. $\Box$
\\

Recall two frames $\fc=\{f_1,\ldots,f_m\}$ and $\gc=\{g_1,\ldots,g_m\}$ for the same Hilbert space $H$ are said {\em equivalent} if there is an invertible
operator $T:H\rightarrow H$ so that $g_k=Tf_k$, for all $1\leq k\leq m$ (see \cite{Bal99,HanLarson00}).The property of being phase retrievable is invariant
among equivalent frames, as the following lemma shows.

\begin{lemma}
\label{lem3}
Assume $\fc=\{f_1,\ldots,f_m\}$ is a phase retrievable frame for $H$. Then
\begin{enumerate}
\item For any invertible operator $T:H\rightarrow H$ and non-zero scalars $z_1,\ldots,z_m\in K$, the frame $\gc=\{g_1,\ldots,g_m\}$ defined by
$g_k=z_k Tf_k$, $1\leq k\leq m$, is also phase retrievable;
\item For any invertible operator $T:H\rightarrow H$, the equivalent frame $\gc=\{g_1,\ldots,g_m\}$ defined by $g_k=Tf_k$, $1\leq k\leq m$ is also
phase retrievable;
\item The canonical dual frame $\tilde{\fc}=\{\tilde{f}_1,\ldots,\tilde{f}_m \}$ is also phase retrievable, where $\tilde{f}_k=S^{-1}f_k$, $1\leq k\leq m$;
\item The associated Parseval frame $\fc^{\#}=\{f_1^{\#},\ldots,f_m^{\#} \}$ is also phase retrievable, where $f_k^{\#}=S^{-1/2}f_k$, $1\leq k\leq m$;
\item Any finite set of vectors $\gc\subset H$ so that $\fc\subset \gc$ is a phase retrievable frame;
\item If $\gc\subset H$ is not a phase retrievable frame then any subset $\hc\subset \gc$ is also not a phase retrievable frame.
\end{enumerate}
\end{lemma}

{\bf Proof of Lemma \ref{lem3}}

(i) Note that each $z_k\neq 0$ and hence $\gc$ is also frame. Let $\alpha_{\gc}:\hat{H}\rightarrow\R^m$ be the nonlinear map associated to $\gc$,
 $(\alpha_{\gc}(x))_k=|\ip{x}{g_k}|^2$. If $x,\in\hat{H}$ are so that $\alpha_{\gc}(x)=\alpha_{\gc}(y)$ then $\nonlin(T^*x)=\nonlin(T^*y)$.
Since $\fc$ is phase retrievable it follows $T^*x=T^*y$ and hence $x=y$. (Note any operator $R:H\rightarrow H$ lifts to a unique
operator $R:\hat{H}\rightarrow\hat{H}$ that is denoted using the same letter).

(ii)-(iv) follows from (i). Claims (v) and (vi) are obvious. $\Box$

\begin{remark}
While the claim (vi) in previous Lemma is obvious, the following question is not so obvious.
Assume $\fc$ is a phase retrievable frame in the real case (that is $H=\R^n$). We know $m\geq 2n-1$. 
Assume $m>2n-1$. Does there always exist a subset $\gc\subset\fc$ so that $\gc$ is a phase retrievable frame?
Interestingly the answer to this question is negative. The following example shows this phenomenon (similar example
was proposed by \cite{Cah12}).
\end{remark}
\begin{example}
Consider $H=\R^3$ and $m=6$ where the 6 vectors are:
\begin{equation}
f_1= \left[ \begin{array}{c} 1 \\ 0 \\ 0 \end{array} \right] ~,~
f_2= \left[ \begin{array}{c} 0 \\ 1 \\ 0 \end{array} \right] ~,~
f_3= \left[ \begin{array}{c} 0 \\ 0 \\ 1 \end{array} \right] ~,~
f_4= \left[ \begin{array}{c} 1 \\ 1 \\ 0 \end{array} \right] ~,~
f_5= \left[ \begin{array}{c} 1 \\ 0 \\ 1 \end{array} \right] ~,~
f_6= \left[ \begin{array}{c} 0 \\ 1 \\ 1 \end{array} \right] 
\end{equation}
The associated rank-1 operators $F_k=f_kf_k^T$, $1\leq k\leq 6$, belong to the linear space of symmetric $3\times 3$ matrices $Sym(\R^3)$.
Note the $Sym(\R^3)$ is a real vectors space of dimension 6. The Gram matrix $G^{(2)}$ associated to $\{F_1,\ldots,F_6\}$ is a $6\times 6$ symmetric matrix
of entries $G^{(2)}_{k,l} = \ip{F_k}{F_l}=|\ip{f_k}{f_l}|^2$, which are the square of the entries of Gram matrix associated to $\fc$. Explicitely $G^{(2)}$ is given by
\begin{equation}
\label{eq:G2}
G^{(2)} =  \left[ \begin{array}{cccccc}
1 & 0 & 0 & 1 & 1 & 0 \\
0 & 1 & 0 & 1 & 0 & 1 \\
0 & 0 & 1 & 0 & 1 & 1 \\
1 & 1 & 0 & 4 & 1 & 1 \\
1 & 0 & 1 & 1 & 4 & 1 \\
0 & 1 & 1 & 1 & 1 & 4
\end{array} \right] 
\end{equation}
Its determinant is $det(G^{(2)})=8$. Hence $\{F_1,\ldots,F_6\}$ is a basis for $Sym(\R^3)$ and thus $ker(\ac)=\{0\}$ which implies $\fc$ is a phase retrievable frame.
On the other hand consider any subset $\gc$ of 5 vectors of $\fc$. It is easy to check $\gc$ is a frame for $\R^3$. Howeverfor each $\gc$ there is a subset of $3$ elements
that is not linearly independent, hence cannot span $\R^3$. This fact together with Corollary 2.6 from \cite{BCE06} proves that $\gc$ is not phase retrievable.  
Thus we constructed a frame $\fc$ of 6 vectors (which is more than the critical cardinal $2n-1=5$) so that any subset is not phase retrievable.
\end{example}
\space{5mm}

We are now ready to present the proof of Theorem \ref{th1}.

{\bf Proof of Theorem \ref{th1} }

Assume $\fc$ is a phase retrievable frame. Then equation (\ref{eq:R2}) is satisfied for some $a_1>0$.
Let $B$ be the upper frame bound for $\fc$. Then set:
\begin{equation}
\label{eq:rho}
\rho = \min(\frac{1}{\sqrt{m}} , \frac{a_1}{4(3B+2)^{3/2}})
\end{equation}

We will find a $\rho>0$ so that (\ref{eq:R2}) is satisfied for any set $\fc'=\{f_1',\ldots,f_m'\}$ with $\norm{f_k-f_k'}<\rho$. 
Let $0<A\leq B<\infty$ be the frame bounds of $\fc$ and let $L'(\xi)$ denote the right hand side in (\ref{eq:R2}) associated to $\fc'$.
We compute

\begin{eqnarray*}
 |\ip{L(\xi)\eta}{\eta} - \ip{L'(\xi)\eta}{\eta}| & \leq & \sum_{k=1}^m |\,|\ip{\Phi_k\xi}{\eta}|^2 - |\ip{\Phi_k'\xi}{\eta}|^2\,|
\leq  \sum_{k=1}^m \left( |\ip{\Phi_k\xi}{\eta}|+|\ip{\Phi_k'\xi}{\eta}| \right) \,|\ip{(\Phi_k-\Phi_k')\xi}{\eta}| \\
 & \leq & \left( \sum_{k=1}^m |\ip{\Phi_k\xi}{\eta}| + \sum_{k=1}^m |\ip{\Phi_k'\xi}{\eta}| \right) \max_{1\leq k\leq m}  |\ip{(\Phi_k-\Phi_k')\xi}{\eta}| 
\end{eqnarray*}
Fix $\xi\in\R^{2n}$ with $\norm{\xi}=1$. Then
\[ \max_{\norm{\eta}=1} \sum_{k=1}^m |\ip{\Phi_k \xi}{\eta}| = \sum_{k=1}^m \ip{\Phi_k \xi}{\frac{\xi}{\norm{\xi}}} \leq B \norm{\xi} \]
Thus for any $\xi,\eta\in\R^{2n}$,
\[ \sum_{k=1}^m |\ip{\Phi_k\xi}{\eta}| \leq B \norm{\xi}\norm{\eta} ~,~\sum_{k=1}^m |\ip{\Phi_k'\xi}{\eta}| \leq B' \norm{\xi}\norm{\eta} \]
where $B'$ is the upper frame bound of $\fc'$. 
On the other hand we bound
\[ |\ip{(\Phi_k-\Phi_k')\xi}{\eta}| \leq \norm{\Phi_k-\Phi_k'}\norm{\xi}\norm{\eta} \]
According to Lemma 3.12 (4) from \cite{Bal13a},
$\norm{\Phi_k-\Phi_k'} = \norm{F_k-F_k'}$, where $F_k=f_kf_k^*$ and $F_k'=f_k'f_k'^*$. Note $F_k-F_k'\in \Soneone(\C^n)$
and $F_k-F_k' = \outp{f_k-f_k'}{f_k+f_k'}$. Now using Lemma 3.7 (4) from \cite{Bal13a}, we obtain
\[ \norm{F_k-F_k'} \leq \norm{F_k-F_k'}_1 = \sqrt{\norm{f_k-f_k'}^2\norm{f_k+f_k'}^2 + (imag(\ip{f_k-f_k'}{f_k+f_k'}))^2} 
\leq \sqrt{2}\norm{f_k-f_k'}\norm{f_k+f_k'} \]
where $\norm{T}_1$ is the nuclear norm (the sum of its singular values) of $T$.
Next notice $\norm{f_k+f_k'}\leq\norm{f_k}+\norm{f_k'}\leq \sqrt{B}+\sqrt{B'}\leq\sqrt{2(B+B')}$. Putting all the estimates together we obtain:
\[  |\ip{L(\xi)\eta}{\eta} - \ip{L'(\xi)\eta}{\eta}| \leq 2(B+B')^{3/2} \left(\max_{1\leq k\leq m}\norm{f_k-f_k'} \right)
\norm{\xi}^2\norm{\eta}^2 \]
Thus
\[ L'(\xi) \geq (a_1 - 2(B+B')^{3/2}\rho)\norm{\xi}^2 1_{\R^{2n}} \]
Finally we obtain an estimate of $B'$ in terms of $B$, $\rho$ and $m$. This estimate can be further refined, but we do not need 
such a refinment  for this proof. Let $\delta_k=f_k'-f_k$. Then
\begin{equation}
\label{eq:Bprime}
 \sum_{k=1}^m |\ip{x}{f_k'}|^2 = \sum_{k=1}^m |\ip{x}{f_k}+\ip{x}{\delta_k}|^2 \leq 2 \left( \sum_{k=1}^m |\ip{x}{f_k}|^2 +
\sum_{k=1}^m|\ip{x}{\delta_k}|^2 \right) = 2(B+m\max_{k}\norm{\delta_k}^2)\norm{x}^2
\end{equation}
Due to (\ref{eq:rho}) we obtain $B'=\sup_{\norm{x}=1}\sum_{k=1}^m |\ip{x}{f_k'}|^2 \leq 2(B+1)$. 
This bound implies that
\[ L'(\xi) \geq \frac{a_1}{2} \norm{\xi}^2 1_{\R^{2n}} \]
and hence $\fc'$ is phase retrievable. ~$\Box$

\section{Critical case $m=4n-4$}
\label{sec3}
This section comments on the recent construction by Bodmann and Hammen \cite{BH13} of a $4n-4$ phase retrievable frame in $\C^n$.
Their construction is as follows. Fix $a\in\R\setminus\pi\QQ$, an irrational multiple of $\pi$. The frame set $\fc$ is given by a union of two sets,
 $\fc=\fc_1\cup\fc_2^a$, where $\fc_1$ constains the following $2n-3$ vectors:
\begin{equation}
\fc_1 =\{f_k^1= \left[ \begin{array}{ccccc}
1 & \mbox{$e^{2\pi i (k+1)/(2n-1)}$} & \mbox{$e^{2\pi i (k+1)2/(2n-1)}$} & \cdots & \mbox{$e^{2\pi i (k+1)(n-1)/(2n-1)}$}
\end{array} \right]^T~~,~~1\leq k\leq 2n-3 \}
\end{equation}
and $\fc_2^a$ contains the following $2n-1$ vectors:
\begin{equation}
\fc_2^a = \{ f_k^2 = \left[ \begin{array}{ccccc}
1 & \mbox{$z_k$} & \mbox{$z_k^2$} & \cdots & \mbox{$z_k^{n-1}$}
\end{array} \right]^T~,~1\leq k\leq 2n-1  \}
\end{equation}
where
\begin{equation}
z_k = \frac{sin\left(\frac{\pi}{2n-1}\right)}{sin(a)}e^{i\frac{k-1}{2n-1}} - e^{i\left(\frac{\pi}{2n-1}-\frac{a}{2}\right)}
\frac{sin\left(\frac{\pi}{2n-1}-\frac{a}{2}\right)}{sin(a)} 
\end{equation}
The proof that $\fc$ is a phase retrievable frame is based on a result by P. Jamming from \cite{Jm10}.
Our Theorem \ref{th1} proves that, in fact, $\fc$ remains phase retrievable for a small perturbation. Since $f_k^2$ depends continuously on $a$,
it follows that the set $\R\setminus\pi\QQ$ can be replaced by a much larger set of real numbers that includes most of rational multiples of $\pi$.
Going through the proof of Theorem 2.3 in \cite{BH13} , and in particular of Lemma 2.2, the only requirement on $a$ is that, any set of $2(n-1)$
complex numbers cannot be simultaneoulsy symmetric with respect to the real line and to a line of angle $a$ passing through the origin.
This phenomenon happens for any $n$ when $a$ is an irrational multiple of $\pi$. However, for a fixed $n$, only finitely many values of
$a$ may allow such a symmetry. In fact when such a symmetric set of $2(n-1)$ complex numbers exists, $a=\pi \frac{p}{q}$ for some $q\leq 2(n-1)$. 
Thus the frame set above $\fc=\fc_1\cup\fc_2^a$ remains phase retrievable for all values of $a$ except a finite set of the form
$\{ \pi\frac{p}{q}~,~0\leq p\leq 2q\leq 4(n-1)\}$.

\section{Non Phase-Retrievable Frames}
\label{sec4}

Consider now the case when $m$ is "small". The conjecture in \cite{BCMN13a} reads that for $m<4n-4$ there is no phase retrievable frame.
In this section we comment on a partial result supporting this conjecture. The main result of this section is the following
\begin{proposition}
\label{prop1}
Fix the $n$-dimensional Hilbert space $H$. Denote by $m^*(n)$ the critical cardinal $m^*(n)=4n-4$.
Assume the following statement holds true for some $m<m^*(n)$:

(O) For any frame set $\fc=\{f_1,\ldots,f_m\}\subset H$ that is not a phase retrievable frame for $H$ there exists a $\rho>0$ so that 
any other set $\fc'=\{f_1',\ldots,f_m'\}$ with $\norm{f_k-f_k'}<\rho$, $1\leq k\leq m$, is also not phase retrievable.

Then any subset $\gc=\{g_1,\ldots,g_m\}\subset H$ of $m$ vectors in $H$ is not a phase retrievable frame.
\end{proposition}

\begin{remark}
Before presenting its proof, we make the following remark. While not proving the fulll conjecture, this result reduces the proof
of the $4n-4$ conjecture to a stability result for non phase retrievable frames. Additionally the result holds true even if the critical cardinal
is not $4n-4$. However the author recognizes this is just a partial result that does not prove the full conjecture. 
\end{remark}

\begin{remark}
Note the set $\fc$ is supposed to be frame. If (O) holds for any set of $m$ vectors of $H$, one can use the trivial set $\{0,0,\ldots,0\}$
of $m$ vectors. Then if (O) holds for this special set, then any set of $m$ vectors whose norms are less than some $\rho>0$ is
non phase retrievable frame. By scaling we obtin immediately that any $m$-set of vectors is not phase retrievable frame.
\end{remark}

{\bf Proof of Proposition \ref{prop1}}

First note that if $m<2n$ any set $\fc$ of $m$ elements cannot be a phase retrievable frame (conform \cite{Fink04}).
Hence we can assume $m\geq 2n$.
Let 
$H^m = H\times H\times\cdots\times H$ denote the $m$-product space endowed with the topology induced by the norm 
$$\norm{\xb}_{H^m} = \max_{1\leq k\leq m} \norm{x^k}~~,~~\xb=(x^1,\ldots,x^m)\in H^m $$
Note that $H^m$ is homeomorphic with  $\C^{nm}$ endowed with the usual Euclidian norm. 
Let $F_n^m$ denote the set of frames for $H$ with $m$ elements. $F^m_n$ is an open set
in $H^m$ since each frame set is stable under small perturbations: for instance this can be seen using an estimate similar to
equation (\ref{eq:Bprime}) that bounds below $\sum_{k=1}^m|\ip{x}{f_k'}|^2$:
\begin{eqnarray*}
 \sum_{k=1}^m |\ip{x}{f_k'}|^2 & \geq & \sum_{k=1}^m |\ip{x}{f_k}|^2 - \left(2 \sqrt{m}\left(\sum_{k=1}^m|\ip{x}{f_k}|^2 \right)^{1/2}
\left(\max_k \norm{f_k-f_k'}\right) \norm{x} + m \left(\max_k \norm{f_k-f_k'}\right)^2 \norm{x}^2 \right) \\
 & \geq & (A-2\sqrt{mB}\rho-m\rho^2)\norm{x}^2
\end{eqnarray*}
Thus for some sufiiciently small $\rho$, $A-2\sqrt{mB}\rho-m\rho^2>0$ and $\{f_1',\ldots,f_m'\}$ is also frame when $\norm{f_k-f_k'}<\rho$, $1\leq k\leq m$.


Assume the hypothesis (O) holds true. 
 Let $N_n^m$ denote the set of non phase-retrievable frames of $m$ vectors in $H$. Thus $N_n^m\subset F^m_n\subset  H^m$
is an open set in $H$ by hypothesis (O).

On the other hand the complement $\Gamma^m_n:=F^m_n\setminus N_n^m$ represents the set of phase retrievable frames. Theorem \ref{th1}
shows $\Gamma^m_n$ is open in $H^m$. 

Now let us show the set of frames $F^m_n$ is connected in $H^m$. Firstly two equivalent frames are connected by path as shown
e.g. in \cite{Bal99}. We will show that any two frames of $m$ elements for the $n$ dimensional Hilbert space $H$ can be
connected by a continuous path (in fact two segments of line), when $m\geq 2n$. Let $\fc_1=\{f_1^1,\ldots,f_m^1\}$ and
$\fc_2=\{f_1^2,\ldots,f_m^2\}$ be two $m$-frames. Let $I=\{k_1,\ldots,k_n\}$ be the $n$-set so that 
$\fc_1[I]=\{f^1_{k_1},\ldots,f^1_{k_n}\}$ is a linearly independent subset of $\fc_1$. Let $J=\{j_1,\ldots,j_{m-n}\}$ be
a $m-n$-set so that $\fc_2[J]=\{ f^2_j~,~j\in J \}$ is a frame for $H$. Let $\gamma:I^c\rightarrow J$ 
and $\delta:J^c\rightarrow I$ be two bijective maps where $I^c=\{1\leq k\leq m\}\setminus I$ and 
$J^c=\{1\leq j\leq m\}\setminus J$ are the complement sets of $I$ and $J$ respectively. We build a piecewise linear path $\beta:[-1,1]\rightarrow H^m$
connecting $\fc_1$ to $\fc_2$ as follows: For $-1\leq t\leq 0$,
\[ (\beta(t))_k = \left\{ \begin{array}{ccl} 
\mbox{$f^1_k$} &  if  & \mbox{$k\in I$} \\
\mbox{$-tf^1_k + (t+1)f^2_{\gamma(k)}$} & if & \mbox{$k\in I^c$}
\end{array} \right. \]
For $0\leq t\leq 1$,
\[ (\beta(t))_k = \left\{ \begin{array}{ccl} 
\mbox{$f^2_j$} &  if  & \mbox{$j\in J$} \\
\mbox{$tf^1_{\delta(j)} + (1-t)f^2_j$} & if & \mbox{$j\in J^c$}
\end{array} \right. \]
One can easilty check that $\beta(t)$ is a frame for each $-1\leq t\leq 1$, and $\beta(-1)=\fc_1$, $\beta(1)=\fc_2$. This proves
the set of frames $F^m_n$ is path connected, hence connected.

We obtained that the connected set $F^m_n$ can be partitioned into two opens sets $\Gamma^m_n$ and $N_n^m$.
It follows that one of the two sets must be the empty set. However we can always construct a non phase retrievable frame, for instance
$\fc=\{e_1,\ldots,e_n,e_n,\ldots,e_n\}$ where $\{e_1,\ldots,e_n\}$ is a basis of $H$ and the vector $e_n$ is repeated
a total of $m-n+1$ times. This shows $\Gamma^m_n$ must be empty. Thus any set $\gc\subset H$ of $m$ vectors cannot
be a phase retrievable frame.~~$\Box$

\section*{ACKNOWLEDGMENTS}
The author has been partially supported by the National Science Foundation under NSF DMS-1109498 grant. 


\end{document}